\documentclass{cedram-aif}
\usepackage{amsmath}  
\usepackage{amssymb, amsthm,thmtools,thm-restate,enumerate,verbatim,bm}
\usepackage{mathrsfs}
\usepackage[usenames,dvipsnames,svgnames,table]{xcolor}
\definecolor{darkgreen}{rgb}{0.0,0.9,0}
\usepackage[pdfauthor=Lyons\ and\ Zhai,bookmarksnumbered,hyperfigures,colorlinks=true,citecolor=BrickRed,linkcolor=BrickRed,urlcolor=BrickRed,pdfstartview=FitH]{hyperref}
\usepackage[inline,nolabel]{showlabels}

\usepackage{microtype}

\providecommand{\norm}[1]{\left\lVert#1\right\rVert}

\equalenv{cor}{coro}
\equalenv{lem}{lemm}
\equalenv{pro}{prop}

\theoremstyle{definition}

\equalenv{defn}{defi}
\equalenv{eg}{exam}

\theoremstyle{remark}

\equalenv{remk}{rema}

\newenvironment{proofof}[1]{{\medbreak\noindent \em Proof of #1.  }}{\hfill\qed\medbreak}

\let\P\relax
\DeclareMathOperator{\P}{\mathbf{P}\mathopen{}}
\DeclareMathOperator{\E}{\mathbf{E}\mathopen{}}

\def\Psub_#1{\P_{\! #1}}

\def\Pbig#1{\P\mkern-.5mu\bigl[#1\bigr]}
\def\Psubbig_#1#2{\Psub_{#1}\mkern-1.5mu\bigl[#2\bigr]}

\def\Psubbigg_#1#2{\Psub_{#1}\mkern-1.5mu\biggl[#2\biggr]}

\def\Ebig#1{\E\mkern-1.5mu\bigl[#1\bigr]}
\def\Esubbig_#1#2{\E_{#1}\mkern-1.5mu\bigl[#2\bigr]}
\def\EsubBig_#1#2{\E_{#1}\mkern-1.5mu\Bigl[#2\Bigr]}
\def\Esupbig^#1#2{\E^{#1}\mkern-1.5mu\bigl[#2\bigr]}
\def\EsupBig^#1#2{\E^{#1}\mkern-1.5mu\Bigl[#2\Bigr]}

\def\Esubbigg_#1#2{\E_{#1}\mkern-1.5mu\biggl[#2\biggr]}
\def\EBig#1{\E\mkern-1.5mu\Bigl[#1\Bigr]}
\def\EBigg#1{\E\mkern-1.5mu\Biggl[#1\Biggr]}

\def\dfn#1{\textit{\textbf{#1}}}

\def\Z{{\mathbb Z}}
\def\R{{\mathbb R}}
\def\C{{\mathbb C}}

\def\cbuldot{{\raise.25ex\hbox{$\scriptscriptstyle\bullet$}}}

\def\bfo{{\bf 1}}

\def\changecomma#1,{#1,\,}
\def\bigchangecomma#1,{#1,\;}
\def\leftchangecomma#1,{#1,\ }

\def\sF{\bm{F}}  
\def\sG{\bm{G}}
\def\sH{\bm{H}}

\def\disk{{\mathbb D}}
\def\T{{\mathbb T}}
\def\dfnterm#1{\dfn{#1}}

\def\bz{{\bf z}}

\def\A#1#2#3{\|a^{(#1)}\|_{#2}^{#3}}

\def\Seq#1{\langle #1 \rangle}
\def\st{\, ; \;}  
\def\doi#1{\url{http://dx.doi.org/#1}}

\let\emptyset\varnothing

\def\thmenv#1#2#3{\begin{#1} \label{#1:#2} #3 \end{#1}}
\def\richthmenv#1#2#3#4{\begin{#1}[#3] \label{#1:#2} #4 \end{#1}}

\def\procl#1.#2 #3\endprocl{%
       \ifx#1t\thmenv{thm}{#2}{#3}\fi
       \ifx#1l\thmenv{lem}{#2}{#3}\fi
       \ifx#1p\thmenv{pro}{#2}{#3}\fi
       \ifx#1c\thmenv{cor}{#2}{#3}\fi
       \ifx#1d\thmenv{dfn}{#2}{#3}\fi
       \ifx#1g\thmenv{conj}{#2}{#3}\fi
       \ifx#1q\thmenv{question}{#2}{#3}\fi
       \ifx#1r\thmenv{remk}{#2}{{\rm #3}}\fi
    }%

\def\rprocl#1.#2 #3 #4\endprocl{%
       \ifx#1t\richthmenv{thm}{#2}{#3}{#4}\fi
       \ifx#1l\richthmenv{lem}{#2}{#3}{#4}\fi
       \ifx#1p\richthmenv{pro}{#2}{#3}{#4}\fi
       \ifx#1c\richthmenv{cor}{#2}{#3}{#4}\fi
       \ifx#1d\richthmenv{dfn}{#2}{#3}{#4}\fi
       \ifx#1g\richthmenv{conj}{#2}{#3}{#4}\fi
       \ifx#1q\richthmenv{question}{#2}{#3}{#4}\fi
       \ifx#1r\richthmenv{remark}{#2}{#3}{{\rm #4}}\fi
    }%

\def\rref#1.#2/{%
      \ifx #1sSection~\ref{s.#2}\fi
      \ifx #1SSubsection~\ref{S.#2}\fi
      \ifx #1tTheorem~\ref{thm:#2}\fi
      \ifx #1lLemma~\ref{lem:#2}\fi
      \ifx #1cCorollary~\ref{cor:#2}\fi
      \ifx #1pProposition~\ref{pro:#2}\fi
      \ifx #1dDefinition~\ref{dfn:#2}\fi
      \ifx #1gConjecture~\ref{conj:#2}\fi
      \ifx #1qQuestion~\ref{question:#2}\fi
      \ifx #1rRemark~\ref{remark:#2}\fi
      \ifx #1aAppendix~\ref{a.#2}\fi
      \ifx #1fFigure~\ref{f.#2}\fi
      \ifx #1e(\ref{e.#2})\fi
      \ifx #1b\citet{#2}\fi
      \ifx #1B\citep{#2}\fi
        }

\newcommand\citep{\cite}
\newcommand\citet{\cite}

\def\refbmulti#1{\citet{#1}}

\def\rlabel #1 #2{\begin{equation} \label{#1} #2 \end{equation}}

\def\rproof{\begin{proof}}

\def\eqaln#1{\begin{align*} #1 \end{align*}}
\def\eqalign#1{\begin{align*} #1 \end{align*}}

\def\beginbulletitems{\begin{itemize}}

\def\endbulletitems{\end{itemize}}

\def\textfrac{\frac}
\newcommand\sfrac[2]{#1/#2}

\def\Qed{\end{proof}}

\def\bsection#1#2{\bigbreak\section{#1}\label{#2}}

\title{Zero Sets for Spaces of Analytic Functions}
\alttitle{Ensembles de z\'eros pour des espaces de fonctions analytiques}

\author{\firstname{Russell} \lastname{Lyons}}
\address{Department of Mathematics\\
831 E. 3rd St.\\
Indiana University\\
Bloomington, IN 47405-7106 (USA)}
\email{rdlyons@indiana.edu}

\author{\firstname{Alex} \lastname{Zhai}}
\address{Department of Mathematics\\
Stanford University\\
450 Serra Mall, Building 380\\
Stanford, CA 94305 (USA)}
\email{azhai@stanford.edu}

\thanks{R.L.\ partially supported by the National
Science Foundation under grant DMS-1612363.
A.Z.\ supported by a Stanford Graduate Fellowship.
Part of this work was done while both authors were visiting
Microsoft Research, Redmond.}


\altkeywords{Bergman, Bargmann, Fock, gaussienne, al\'eatoire}


\begin{document}

\begin{abstract}
We show that under mild conditions, a Gaussian analytic function $\bm F$
that a.s.\ does not belong to a given weighted Bergman space or
Bargmann--Fock space has the property that a.s.\ no non-zero function in that space
vanishes where $\bm F$ does. This
establishes a conjecture of Shapiro (1979) on Bergman spaces
and allows us to resolve a question of Zhu (1993) on Bargmann--Fock spaces.
We also give a similar result on the union of two (or more) such zero sets, thereby
establishing another conjecture of Shapiro (1979) on Bergman spaces and allowing us to
strengthen a result of Zhu (1993) on Bargmann--Fock spaces.
\end{abstract}

\begin{altabstract}
On montre que sous des conditions faibles, une fonction analytique
gaussienne $\bm F$ qui n'appartient pas p.s.\ \`a un espace pond\'er\'e
de Bergman ou de Bargmann--Fock donn\'e a p.s.\ la propri\'et\'e qu'il
n'existe pas de fonction non-nulle dans cette espace qui s'annule o\`u
$\bm F$ s'annule.
Ceci d\'emontre une conjecture de Shapiro (1979) sur les espaces de Bergman
et nous permet de r\'esoudre une question de Zhu (1993) sur les espaces de
Bargmann--Fock.
On donne aussi un r\'esultat similaire sur la r\'eunion de deux (ou plus) tels
ensembles de z\'eros, montrant ainsi une autre conjecture de Shapiro
(1979) sur les espaces de Bergman et nous permettant de renforcer un
r\'esultat de Zhu (1993) sur les espaces de Bargmann--Fock.
\end{altabstract}

\maketitle

\bsection{Introduction}{s.intro}

Zeros of Gaussian analytic functions were originally studied by
Paley and Wiener \rref b.PW:FT/, Kac \refbmulti{Kac43,Kac43corr}, and Rice
\refbmulti{Rice44,Rice45}.
Since then, many more mathematicians and physicists have been interested in
such zero sets.
For some of the history, see \rref b.Sodin:ECM/ and \rref
b.HKPV:book/. Those sources also give surveys of certain aspects of zero sets of
Gaussian analytic functions as random objects.
The topic of the present paper, however, is not mainly zero sets of
Gaussian analytic functions as random objects, but as tools to understand zero sets
in standard spaces of analytic functions.
In particular, we consider the (weighted) Bergman spaces in the unit disk and the (weighted)
Bargmann--Fock spaces in the entire plane, for which we give a unified treatment.
In \autoref{sec:hist}, we give a brief history of what is known for zero sets of functions in these spaces, focused on results relevant to ours. More can be found in
Chapter 4 of \rref b.HKZ:book/ and
Chapter 4 of \rref b.Duren:book/, which are devoted to zero sets of Bergman
spaces, and
Chapter 5 of \rref b.Zhu:book/, which is devoted to zero sets of Bargmann--Fock
spaces.

Let $\mu$ be a finite measure on $(0, \infty)$, not identically 0.
Write $r_\mu := \inf \{ r \st \mu(r, \infty) = 0\} \in (0, \infty]$, and assume that $\mu\bigl(\{r_\mu\}\bigr) = 0$. For $p \in (0, \infty)$, write $A^p(\mu)$ for the set of
analytic functions $f$ defined for $|z| < r_\mu$ that
satisfy
\[
\int_0^{r_\mu} \int_0^1 |f(r e^{2\pi i \theta})|^p \,d\theta
\,d\mu(r) < \infty
\,.
\]
When $r_\mu = 1$, these spaces are referred to as \dfn{weighted Bergman
spaces}, whereas when $r_\mu = \infty$, they are called \dfn{weighted
Bargmann--Fock spaces}. Clearly all spaces $A^p(\mu)$ when $r_\mu$ is finite
are isomorphic to weighted Bergman spaces.
Denote the unit disk by $\disk := \{z \st |z| < 1\}$.
The unweighted Bergman spaces $A^p(\disk)$
correspond to $d\mu(r) = 2r \bfo_{[0, 1]}(r) \,dr$.
The most-studied weights are $d\mu(r) = 2(1-r^2)^\alpha \bfo_{[0, 1]}(r)
\,dr$ ($\alpha > -1$),
in which case the corresponding Bergman spaces are denoted
$A^p_\alpha(\disk)$. By contrast, the most-studied
Bargmann--Fock spaces are defined
differently, with $\mu$ depending on $p$, namely,
$d\mu(r) = p\alpha r\, e^{- p \alpha r^2/2}\,dr$ ($\alpha > 0$),
in which case the corresponding Bargmann--Fock
spaces are denoted $B^p_\alpha(\C)$. An older definition of $B^p_\alpha(\C)$ was used by \rref
b.Zhu:zerosFock/, where $\mu$ did not depend on $p$; in our notation, this was the space
$B^p_{2\alpha/p}(\C)$. By \cite[Theorem 2.10]{Zhu:book}, $B^p_\alpha(\C) \subsetneq B^q_\alpha(\C)$
for $0 < p < q < \infty$.

A \dfn{standard complex Gaussian} random variable is one whose density
with respect to Lebesgue measure on $\C$ is $z \mapsto e^{-|z|^2}/\pi$.
We always consider the zero set $Z(f)$ of an analytic function $f$ as a
multiset or a sequence, where each zero $w$ is listed with its
multiplicity, which is $m$ if $z \mapsto f(z)/(z-w)^m$ is analytic and does not vanish at $w$.

Our main result is the following.

\procl t.main
Let $\mu$ be a finite measure on $(0, \infty)$ with $\mu\bigl(\{r_\mu\}\bigr) = 0$. Let $p \in (0, \infty)$.
Suppose that $a_n \ge 0$ satisfy $\limsup_{n \to\infty} a_n^{1/n} < \infty$
and $r \mapsto \sum_{n = 0}^\infty a_n^2 r^{2n} \notin L^{p/2}(\mu)$.
Let $\sF(z) := \sum_{n=0}^\infty a_n \zeta_n z^n$ for $|z| < r_\mu$, where
$\zeta_n$ are independent complex Gaussian random variables.
Then a.s.\ the only analytic function $f \in A^p(\mu)$ with $Z(f) \supseteq
Z(\sF)$ is $f \equiv 0$.
\endprocl

Note that if $r \mapsto \sum_{n = 0}^\infty a_n^2 r^{2n} \notin L^{p/2}(\mu)$ for all $p > p_0$, then by considering a countable set of $p > p_0$, we may conclude that
a.s.\ for all $p > p_0$, the only analytic function $f \in A^p(\mu)$ with $Z(f) \supseteq
Z(\sF)$ is $f \equiv 0$.

The following corollary,
in the special case where $\mu_1 = \mu_2$,
was known for Bergman spaces \citep{horo:0sberg}.
The corollary follows from \rref t.main/ and \rref e.fernAp/.

\procl c.main
Let $R \in (0, \infty]$.
Let $\mu_i$ ($i = 1, 2$) be finite measures with $r_{\mu_i} = R$.
Let $p_i \in (0, \infty)$ ($i = 1, 2$).
Suppose that there exist
$a_n \ge 0$ that satisfy $\limsup_{n \to\infty} a_n^{1/n} < \infty$
and $r \mapsto \sum_{n = 0}^\infty a_n^2 r^{2n} \in L^{p_1/2}(\mu_1)
\setminus L^{p_2/2}(\mu_2)$.
Then there is a function $f \in A^{p_1}(\mu_1)$ such that the only $g \in
A^{p_2}(\mu_2)$ with $Z(g) \supseteq Z(f)$ is
$g \equiv 0$.
\endprocl

We will actually prove a quantitative version of \rref t.main/.
Write $A(\mu)$ for the set of functions that are analytic in $\{z \st |z| <
r_\mu\}$.
For $s \le r_\mu$ and $f \in A(\mu)$, write
$Z_s(f)$ for the multiset of $z$ with $f(z) = 0$ and $0 < |z| < s$.
Denote
\[
\norm{f}_{A^p(\mu, s)}
:=
\Bigl(
\int_0^{s} \int_0^1 |f(r e^{2\pi i \theta})|^p \,d\theta
\,d\mu(r) \Bigr)^{1/p}
.
\]
We also abbreviate
\[
\norm{f}_{A^p(\mu)}
:=
\norm{f}_{A^p(\mu, r_\mu)}
\,.
\]
Given a sequence $\Seq{a_n \st n \ge 0}$,
write $a^{(r)}$ for the sequence $\Seq{a_n r^n \st n \ge 0}$ and $\A r2{}$ for
its $\ell^2$-norm.

\procl t.quant
Let $a_n \ge 0$ satisfy $R^{-1} := \limsup_{n \to\infty} a_n^{1/n} < \infty$ and $a_0 \ne 0$.
Let $\sF(z) := \sum_{n=0}^\infty a_n \zeta_n z^n$ for $|z| < R$, where
$\zeta_n$ are independent complex Gaussian random variables.
Then for all finite measures $\mu$ with $r_{\mu} = R$ and $\mu\bigl(\{R\}\bigr) = 0$, all $p \in (0, \infty)$,
and all $s \in (0, R]$,
\begin{multline} \label{e.quant}
\EBig{\max \Bigl\{ \frac{|f(0)| }{ \norm{f}_{A^p(\mu, s)}} \st 0 \not\equiv
f \in A(\mu),\, Z(f) \supset Z_s(\sF) \Bigr\}}
\\ \le
\frac{\sqrt\pi\, a_0 }{ \bigl(\int_0^s \A {r}2p
\,d\mu(r)\bigr)^{1/p}}
\,.
\end{multline}
\endprocl

\begin{proofof}{\rref t.main/ from \rref t.quant/}
Consider $0 \not\equiv f \in A(\mu)$ with $Z(\sF) \subset Z(f)$; we will
show that $\|f\|_{A^p(\mu)} = \infty$.

Without loss of generality, we may shift the indices of the $a_n$ so that $a_0 \ne 0$, since this does not affect the condition $\|a^{(r)}\|_2 \not\in L^p(\mu)$, and it does not change $Z_R(\sF)$.  Thus, \rref t.quant/ applies.

If $f(0) \ne 0$, then the result follows directly from \rref e.quant/ by taking $s = R$. Otherwise, we may reduce to this case: Let $m$ denote the order of vanishing of $f$ at $0$, and let $g(z) := f(z)/z^m$, so that $g \in A(\mu)$ and $g(0) \ne 0$. We then have $Z(g) \supset Z_R(\sF)$, from which we conclude that $\norm{g}_{A^p(\mu)} = \infty$.
This clearly implies that also
$\norm{f}_{A^p(\mu)} = \infty$, as desired.
\end{proofof}

We also establish the following theorem, which relates to unions of zero sets in the special case $b \equiv -1$ upon observing that $Z(\sF^N -1) = \bigcup_{k=0}^{N-1} Z(\sF - e^{2\pi i k/N})$.

\procl t.main2
Let $\mu$ be a finite measure on $(0, \infty)$ with $\mu\bigl(\{r_\mu\}\bigr) = 0$. Let $p \in (0, \infty)$.
Suppose that $a_n \ge 0$ satisfy $\limsup_{n \to\infty} a_n^{1/n} < \infty$
and $r \mapsto \sum_{n = 0}^\infty a_n^2 r^{2n} \notin L^{p/2}(\mu)$.
Let $\sF(z) := \sum_{n=0}^\infty a_n \zeta_n z^n$ for $|z| < r_\mu$, where
$\zeta_n$ are independent complex Gaussian random variables.
Let $b \in A(\mu)$ and $N$ be a positive integer.
Then a.s.\ the only analytic function $f \in A^{p/N}(\mu)$ with $Z(f) \supseteq
Z(\sF^N + b)$ is $f \equiv 0$.
\endprocl

This establishes the full conjecture of Shapiro
\rref b.Shapiro/ and enlarges the set
of $b$ to which it applies.
Since $Z(\sF\pm1)$
are a.s.\ simple (see \citet[Lemma 28]{PeresVirag}) and $Z(\sF\pm 1)$ are disjoint, we obtain the following corollary.

\procl c.main2
Let $R \in (0, \infty]$.
Let $\mu_i$ ($i = 1, 2$) be finite measures with $r_{\mu_i} = R$ and $\mu_i\bigl(\{R\}\bigr) = 0$.
Let $p_i \in (0, \infty)$ ($i = 1, 2$).
Suppose that there exist
$a_n \ge 0$ that satisfy $\limsup_{n \to\infty} a_n^{1/n} < \infty$
and $r \mapsto \sum_{n = 0}^\infty a_n^2 r^{2n} \in L^{p_1/2}(\mu_1)
\setminus L^{p_2/2}(\mu_2)$.
Then there are functions $f_1, f_2 \in A^{p_1}(\mu_1)$ such that $Z(f_1)
\cap Z(f_2) = \emptyset$ and the only $g \in
A^{p_2/2}(\mu_2)$ with $Z(g) \supseteq Z(f_1) \cup Z(f_2)$ is
$g \equiv 0$.
\endprocl

Again, we prove a quantitative version of \rref t.main2/:

\procl t.quant2
Let $a_n \ge 0$ satisfy $R^{-1} := \limsup_{n \to\infty} a_n^{1/n} < \infty$.
Let $\sF(z) := \sum_{n=0}^\infty a_n \zeta_n z^n$ for $|z| < R$, where
$\zeta_n$ are independent complex Gaussian random variables.
Then for all finite measures $\mu$ with $r_{\mu} = R$ and $\mu\bigl(\{R\}\bigr) = 0$, all $p \in (0, \infty)$,
all $b \in A(\mu)$, all positive integers $N$,
and all $s \in (0, R]$,
\begin{multline*}
\EBig{\max \Bigl\{ \frac{|f(0)|^{1/N} }{ \norm{f}_{A^p(\mu, s)}^{1/N}}
\st 0 \not\equiv
f \in A(\mu),\, Z(f) \supset Z_s(\sF^N + b) \Bigr\}}
\\ \le
\frac{c }{ \bigl(\int_0^s \A {r}2p
\,d\mu(r)\bigr)^{1/p}}
\,,
\end{multline*}
where
\[
c :=
\Bigl(a_0^{4N} (2N)! + 4|b(0)|^2 a_0^{2N} N! + |b(0)|^4\Bigr)^{1/4N}
\Gamma\Bigl(\frac{2N-1}{4N-1}\Bigr)^{\frac{4N - 1}{4N}}
\,.
\]
Therefore, if $\A r2{} \notin L^p(\mu)$, then a.s.\ every $f \in A(\mu)$
with $Z(f) \supset Z(\sF^N + b)$ satisfies $\norm{f}_{A^p(\mu)} = \infty$.
\endprocl

Of course, what allows Gaussian series to have these properties is that
such series have many zeros. A quantitative form of this property is what
lies behind our results.
Recall that by the arithmetic mean-geometric mean inequality (or Jensen's
inequality) and Jensen's
formula, every $f \in A(\mu)$ with $f(0) \ne 0$ satisfies
\begin{align}
\label {e.AMGMJensen}
\begin{split}
\norm{f}_{A^p(\mu)}^p
&=
\int_0^{r_\mu} \int_0^1 |f(r e^{2\pi i \theta})|^p \,d\theta \,d\mu(r)
\\ &\ge
\int_0^{r_\mu} \exp \int_0^1 \log |f(r e^{2\pi i \theta})|^p \,d\theta
\,d\mu(r)
\\ &=
\int_0^{r_\mu} |f(0)|^p \prod_{z \in Z(f)} \max\bigl\{\textfrac{r^p}{|z|^p},
1\bigr\} \,d\mu(r)\,.
\end{split}
\end{align}
In general, this inequality can be very far from an equality; for two
simple examples, consider $f(z) := (1-z)^{-1}$ and $p \ge 2$ or $f(z) :=
e^{(1-z)^{-1}}$ and all $p$.
What we will show, in contrast, is that for $f = \sF$, a.s.\ finiteness of
the right-hand side of \rref e.AMGMJensen/ implies a.s.\ finiteness of the
left-hand side and even finiteness of the expectation of the left-hand side.
This is reminiscent of Fernique's theorem (Appendix \ref{s.fernique}), but the
functional on the right-hand side does not satisfy the hypotheses of
Fernique's theorem.
Moreover, Fernique's theorem gives finiteness of a moment defined in terms
of the original functional, whereas here, the $A^p(\mu)$-norm is, as we
just illustrated, not in any way a function of the right-hand side.

\procl t.quant3
Let $a_n \ge 0$ satisfy $R^{-1} := \limsup_{n \to\infty} a_n^{1/n} < \infty$ and $a_0 \ne 0$.
Let $\sF(z) := \sum_{n=0}^\infty a_n \zeta_n z^n$ for $|z| < R$, where
$\zeta_n$ are independent complex Gaussian random variables.
Then for all finite measures $\mu$ with $r_{\mu} = R$ and
$\mu\bigl(\{R\}\bigr) = 0$ and all $p \in (0, \infty)$, the following are
equivalent:
\begin{enumerate}[\rm (i)]
\item $\int_0^{R} \exp \int_0^1 \log |\sF(r e^{2\pi i \theta})|^p \,d\theta
\,d\mu(r) < \infty$ a.s.;\vadjust{\kern2pt}%
\item $\Ebig{\norm{\sF}_{A^p(\mu)}^p} < \infty$;\vadjust{\kern2pt}%
\item $\Ebig{\int_0^{R} \exp \int_0^1 \log |\sF(r e^{2\pi i \theta})|^p \,d\theta
\,d\mu(r)} < \infty$;\vadjust{\kern3pt}%
\item $\int_0^{R} \exp \int_0^1 \log |\sF(r e^{2\pi i \theta})|^p \,d\theta
\,d\mu(r) < \infty$ with positive probability.
\end{enumerate}
Moreover, for
all $s \in (0, R]$,
\begin{multline} \label{e.quant3}
\EBigg{\frac{|\sF(0)| }{ \Bigl(\int_0^{s} \exp \int_0^1 \log |\sF(r e^{2\pi i \theta})|^p \,d\theta
\,d\mu(r)\Bigr)^{1/p}}}
\\ \le
\frac{\sqrt\pi\, \Gamma(1+p/2)^{1/p}\, a_0 }{ \Ebig{\norm{\sF}_{A^p(\mu, s)}^p}^{1/p}}
\,.
\end{multline}
\endprocl

The equivalence shown here may be surprising; indeed, in discussing his
conjecture, Shapiro
\rref b.Shapiro/ wrote that the arithmetic mean-geometric mean
inequality ``seems to give away too much."

\subsection{History of Zero Sets}\label{sec:hist}

Given a collection $A$ of analytic functions, say that $Z$ is an \dfn{$A$-zero set} if there is some function in $A$ whose zero set equals $Z$.
There is no geometric characterization known for a set of points in $\disk$
to be an $A^p(\disk)$-zero set, but there are necessary
conditions known that are not far from known sufficient conditions. It is
also known that no condition depending solely on the moduli of the points
can be both necessary and sufficient.
For further discussion,
let $\bz$ be a countable multiset in $\disk$ and write
$$
\varphi_{\bz}(r) := \sum_{\substack{z \in \bz, \\ |z| \le r}} \bigl(1 - |z|\bigr)
\,.
$$
The situation for zeros of Bergman functions
contrasts strongly with that for the Hardy spaces,
$$
H^p(\disk)
:=
\{ f \in A^0(\disk) \st \sup_{r < 1} \int_0^1 |f(r e^{2\pi i \theta})|^p
\,d\theta
< \infty \}
\,,
$$
where for all $p \in (0, \infty]$, the \dfn{Blaschke condition}
$$
\varphi_{\bz}(1) < \infty
$$
is necessary and sufficient to be an $H^p(\disk)$-zero set.
For every $p \in (0, \infty)$,
the Blaschke condition is sufficient to be an $A^p(\disk)$-zero set (since
$H^p(\disk) \subset A^p(\disk)$), while
the condition
$$
\sum_{z \in \bz \setminus \{0\}}
\frac{\bigl(1 - |z|\bigr) }{ \log^{1+\epsilon}\bigl(1 -
|z|\bigr)^{-1}} < \infty
$$
is known to be necessary for every $\epsilon > 0$ but not for $\epsilon =
0$ \rref B.horo:0sberg/.
On the other hand,
if a subset of $\bz$ lies on a line (or in a Stolz angle),
then the Blaschke condition for that
subset is also necessary for $\bz$ to be an $A^p(\disk)$-zero set \rref B.SS/.
Combining the preceding results, we deduce that the moduli alone do
not determine whether a point set is an $A^p(\disk)$-zero set.

A set $W$ is called a \dfnterm{set of uniqueness for $A(\mu)$} if
the only $f \in A(\mu)$ with $Z(f) \supseteq W$ is $f \equiv 0$.
Horowitz
\rref b.horo:0sberg/ showed that for $0 < p < q < \infty$, there exists an
$A^p(\disk)$-zero set that is an
$A^q(\disk)$-uniqueness set.
In fact, he showed that if $f \in A^q(\disk)$ with zero set $\Seq{z_k \st k
\ge 1}$
ordered so that $|z_k|$ is increasing and $f(0) \ne 0$, then
\rlabel e.necAp
{
\sup_n n^{-1/q} \prod_{k=1}^n \frac{1 }{ |z_k|} < \infty
\,,
}
whereas for every $p < q$, there is some
$f \in A^p(\disk)$ with zero set $\Seq{z_k \st k \ge 1}$
ordered so that $|z_k|$ is increasing and $f(0) \ne 0$ satisfying
$$
\sup_n n^{-1/q} \prod_{k=1}^n \frac{1 }{ |z_k|} = \infty
\,.
$$
(Since \rref e.necAp/ depends only on the moduli, it is not sufficient to be
a zero set.)
This distinction among the zero sets for different $p$
was refined by Shapiro \rref b.Shapiro/:
for $0 < p < \infty$, there exists $f
\in A^p(\disk)$ whose zeros are not the zeros of any function in
$A^{p+}(\disk)$, where
$A^{p+}(\disk) := \bigcup_{q > p} A^q(\disk)$.\footnote{In that paper, \rref b.Shapiro:wtd/ is cited for the first proof of this existence. However, he seems to have misinterpreted the order of quantifiers.
Instead, the novelty of \rref b.Shapiro:wtd/ was to
extend the allowed set of weights from those in \rref b.horo:0sberg/.}
Shapiro \rref b.Shapiro/ did this by using random (Gaussian) series, as we detail
soon.\footnote{Actually, there was a gap in his proof: in the middle of page 168 where the
quantity $I(r)$ is being
bounded below, going from the integral over $\T$ to $E^\omega(n)$ throws
away a part that may be negative, so the inequality does not follow. Thus,
it seems that our proof of \rref c.main/
is the first valid proof of \citep[Theorem 1 (i) implies (iii)]{Shapiro}.}

Later works
\cite{leblanc,bomash,NW}
considered random angles for fixed moduli,
culminating in the following result.

\procl t.NW
Let $0 < p < \infty$ and $\bz = \Seq{z_n \st n \ge 1}
\subset \disk$. Let $\theta_n$ be independent
uniform $[0, 1]$ random variables. If there exists $\epsilon > 0$ such that
\rlabel e.bomash
{
\int_0^1 e^{p \varphi_{\bz}(r)} \log^{(1+\epsilon)}(1-r)^{-1} \,dr < \infty
\,,
}
then a.s.\ $\Seq{z_n e^{2\pi i \theta_n} \st n \ge 1}$ is an $A^p(\disk)$-zero set.
If $q > p$, then the condition \rref e.bomash/ is not sufficient for
$\Seq{z_n e^{2\pi i \theta_n} \st n \ge 1}$ to be a.s.\ an $A^q(\disk)$-zero set.
\endprocl

The Blaschke condition shows that the union of two $H^p(\disk)$-zero sets
is again an $H^p(\disk)$-zero set.
Horowitz \rref b.horo:0sberg/ also showed that although the union of two
$A^p(\disk)$-zero sets is an $A^{p/2}(\disk)$-zero set (trivially: just
multiply the functions), it need not be an $A^q(\disk)$-zero set if $q >
p/2$.
This was again strengthened by Shapiro
\rref b.Shapiro/ to show that it need not be
an $A^{(p/2)+}(\disk)$-zero set.\footnote{The same gap as noted in the
previous footnote applies to this result, but is filled by the proof of our \rref
c.main2/.}

Many of the above results were extended to weighted Bergman spaces.
For example,
for $(p, \alpha) \in (0, \infty) \times (-1, \infty)$,
Horowitz \rref b.horo:0sberg/
studied the zero sets of the spaces
$A^p_\alpha(\disk)$, showing that they were distinct classes of sets
for pairs with distinct values of $(\alpha+1)/p$, provided that $\alpha \geq 0$.
He asked whether it sufficed that the pairs $(p, \alpha)$ be distinct. The
proviso that $\alpha \geq 0$ was removed by Sedletski\u\i\ \rref b.Sed/.
The full question was answered affirmatively by Sevast{$'$}yanov and
Dolgoborodov \rref b.SD:wtd/.
Our \rref c.main/ easily establishes the
full result of \rref b.SD:wtd/ by using Theorem 1 of \rref b.MM/, which
implies that for $\alpha > -1$, $\,q > 0$, and $c_k \ge 0$,
\[
\int_0^1 \biggl(\sum_{k=0}^\infty c_k r^k\biggr)^q (1 - r)^\alpha \,dr <
\infty
\quad
\text{iff}
\quad
\sum_{n=0}^\infty 2^{-n(\alpha+1)} \biggl(\sum_{k=2^n}^{2^{n+1}-1}
c_k\biggr)^q < \infty\,.
\]
 In the above expressions, we take $c_{2k} = a_k^2$ and $c_{2k+1}
  = 0$, and we make a judicious choice of $a_k$ so that the
  convergence behavior is different for distinct pairs $(p, \alpha)$
  and $(p', \alpha')$. The most interesting case is when $(\alpha +
  1)/p = (\alpha' + 1)/p'$ and $p < p'$, which can be handled by
  taking $q = p/2$, $a_{2^n}^2
  = 2^{2n(\alpha + 1)/p} \cdot n^{-2/p}$, and $a_k = 0$ when $k$ is not
  a power of $2$.

Very little is known about the zero sets of functions in the Bargmann--Fock
spaces, even for $p = 2$.
Zhu \rref b.Zhu:zerosFock/ showed
that if $f \in B^p_\alpha(\C)$ with $f(0) \ne 0$ and we write
$Z(f) = \bz$ as a sequence in increasing order of modulus, then $\inf_n
|z_n|/\sqrt n > 0$.
On the other hand, classical results show that
if $\bz$ satisfies $\sum_n |z_n|^{-2} < \infty$, then
there is some $f \in B^p_\alpha(\C)$ with $Z(f) = \bz$ (see \citet[Theorem
5.3]{Zhu:book}).

The paper
\rref b.CLP/ considered particular stationary
random point processes and showed that for $p = 2$, the critical density
for being a $B^2_1(\C)$-zero set is 1.
Zhu \cite[p.~203]{Zhu:book} gives examples showing that a $B^2_\alpha(\C)$-zero set and a
$B^2_\alpha(\C)$-uniqueness set can differ by just one point, and that for all $p, q \in (0, \infty)$
and $n \in (2/q, \infty) \cap \Z$, for every nontrivial
$B^p_\alpha(\C)$-zero set $W$, removing any $n$ points from $W$ yields a
$B^q_\alpha(\C)$-zero set.

Our results give new proofs of results of Zhu \cite{Zhu:zerosFock,Zhu:book}
and answer his question \cite[pp.~202, 209]{Zhu:book}, showing that the zero sets of
$B^p_\alpha(\C)$ depend on $p$ for fixed $\alpha$; he had shown that they
differ for differing $\alpha$, whether or not $p$ is fixed \cite[Theorem 5.8]{Zhu:book}.
To apply \rref c.main/ to Zhu's question, we use the following result of Stokes \rref b.Stokes/:
for $b \ge 0$,
\rlabel e.stokes
{\lim_{t \to\infty} e^{-t} t^{b} \sum_{n=0}^\infty \frac{t^n}{\Gamma(n+b)}
= 1\,.}
Given $\alpha$ and $p$, set $a_n := \sqrt{\alpha^n/\Gamma(n+2/p)}$ and $\sF(z) := \sum_{n=0}^\infty a_n \zeta_n z^n$,
where $\zeta_n$ are independent complex Gaussian random variables. We
apply \rref e.stokes/ with $b = \sfrac{2}{p}$ and $t = \alpha
r^2$, so that for $q > 0$ and as $r \rightarrow \infty$,
  \[ r e^{-\frac{q \alpha r^2}{2}} \biggl( \sum_{n = 0}^\infty a_n^2 r^{2n}
  \biggr)^{\frac{q}2} \asymp r e^{-\frac{q \alpha r^2}{2}} \bigl(
  e^{\alpha r^2} r^{-\frac{4}{p}} \bigr)^{\frac{q}2} = r^{1 - \frac{2q}{p}}. \]
Then by \rref e.fernAp/ and \rref t.main/, it follows that
a.s.\ $\sF \in \bigcap_{q > p} B^q_\alpha(\C)$ and
$Z(\sF)$ is a $B^p_\alpha(\C)$-uniqueness set.

Similarly, for $b \ge 0$ and $c \in \R$, we have the asymptotic
\[
\lim_{t \to\infty} e^{-t} t^{b} (\log t)^{c} \sum_{n=0}^\infty
\frac{t^n}{\Gamma(n+b) (\log n)^c}
= 1\,.
\]
With $a_n := \sqrt{\alpha^n/\bigl(\Gamma(n+2/p) (\log
  n)^{4/p}\bigr)}$, $\>b = \sfrac{2}{p}$, $\>c = \sfrac{4}{p}$, and $t = \alpha r^2$, we
  obtain by a similar calculation that a.s.\ $\sF \in B^p_\alpha(\C)$
and the only $f \in \bigcup_{q < p} B^q_\alpha(\C)$ with $Z(f)
\supseteq Z(\sF)$ is $f \equiv 0$.

We also strengthen Zhu's result (\rref b.Zhu:zerosFock/ or \cite[Theorem 5.4]{Zhu:book}) that there is a union of two disjoint
$B^p_\alpha(\C)$-zero sets that is a $B^p_\alpha(\C)$-uniqueness
set. Indeed, by \rref t.main2/ and \rref p.flexible/, we can find
disjoint $B^p_\alpha(\C)$-zero sets $Z_1$, $Z_2$ such that the only $f
\in \bigcup_{q > p} B^{q/2}_{2p\alpha/q}(\C)$ with $Z(f) \supseteq Z_1
\cup Z_2$ is $f \equiv 0$; taking $q := 2p$ gives Zhu's result. (Note
that $B^{q/2}_{2p\alpha/q}(\C)$ decreases in $q$ by \cite[Corollary
  2.8]{Zhu:book}.)

\subsection{Shapiro's Approach}\label{sec:shap-app}

Consider $\sF(z) := \sum_{n=0}^\infty a_n \zeta_n z^n$, where $\zeta_n$ are
IID standard complex Gaussian random variables and $a_n > 0$ satisfy
$\limsup_{n \to\infty} a_n^{1/n} \le 1$.
Because $\Ebig{\log^+ |\zeta_0|} < \infty$, we also have
$\limsup_{n \to\infty} |\zeta_n|^{1/n} = 1$ a.s.\ by the Borel--Cantelli
lemma, whence a.s.\ $\sF(z)$ converges for all $z \in \disk$ to an analytic
function.

Let $\mu$ be a finite measure with $r_\mu = 1$.
Write $L^{p+}(\mu) :=
\bigcup_{q > p} L^q(\mu)$.
Shapiro \rref b.Shapiro/ showed that
the following are equivalent:
\begin{enumerate}[(1)]
\item $\bigl(r \mapsto \A r2{} \bigr) \in L^p(\mu) \setminus L^{p+}(\mu)$;
\item a.s.\ $\sF \in A^p(\mu) \setminus A^{p+}(\mu)$;
\item a.s.\ $\sF \in A^p(\mu)$ and the only function in
$A^{p+}(\mu)$ that vanishes everywhere that $\sF$ does is the 0 function.
\end{enumerate}
In addition, he showed that when (1) holds,
\eqalign{
\text{a.s. } &\sF \pm 1 \in A^p(\mu) \text{ and the only function in }
A^{(p/2)+}(\mu)
\\ &\text{ that vanishes on } Z(\sF^2 - 1) \text{ is the 0
function.}
}
He conjectured that the following strengthening holds:
\eqalign{
\bigl(r \mapsto \A r2{} \bigr) &\notin L^p(\mu)
\Longrightarrow  \\ &\hbox{ a.s.\ the only function in $A^{p}(\mu)$
that vanishes on $Z(\sF)$}
\\ &\text{ is the 0 function and the only function in } A^{p/2}(\mu)
\\ &\text{ that vanishes on }
Z(\sF^2 - 1) \text{ is the 0 function.}
}
More generally, he conjectured \rref t.main2/ when $r_\mu = 1$ and $b$
satisfies a certain restriction.

The equivalence of (1) and (2) follows from the following equivalence:
\rlabel e.fernAp
{\bigl(r \mapsto \A r2{} \bigr) \in L^p(\mu)
\quad\Longleftrightarrow \quad
\text{a.s. } \sF \in A^p(\mu)
\,.
}
(This equivalence is valid for $r_\mu = \infty$ as well.)
To see this, note that for each $z \in \disk$, the random variable $\sF(z)$
has the same distribution as $\A {|z|}2{} \zeta_0$.
Thus, Tonelli's theorem yields
\rlabel e.tonelli
{\Ebig{\norm{\sF}_{A^p(\mu)}^p} = \Ebig{|\zeta_0|^p} \cdot \int_0^1 \A r2p \,d\mu(r)
\,.}
The forward implication of \rref e.fernAp/ is now immediate.
The reverse implication is a consequence of \rref e.tonelli/ and Fernique's theorem, which tells
us that if $\sF$ a.s.\ belongs to $A^p(\mu)$, then there
exist some $c_0, c_1 > 0$ such that $\Ebig{\exp\{c_0
\norm{\sF}_{A^p(\mu)}^{c_1}\}} < \infty$.
(See Appendix \ref{s.fernique} for a statement and proof of a general form of Fernique's theorem.)

The usefulness of Shapiro's approach comes partly from his implicit observation
that given $\mu$ and $p$,
there exists
$\bigl(r \mapsto \sum_{n = 0}^\infty a_n^2 r^{2n} \bigr) \in L^{p/2}(\mu)
\setminus \bigcup_{q > p} L^{q/2}(\mu)$. This follows from the lemma in Section 3 of \rref b.Shapiro:wtd/, where he considers analytic functions, not just real power series.
For completeness, we give a short proof and extension here.

\procl p.flexible
Let $\mu$ be a finite measure with $0 < r_\mu \le \infty$; if $r_\mu = \infty$, then assume that
$\int_0^{r_\mu} r^n \,d\mu(r) < \infty$ for every $n \ge 0$.
For all $p \in (0, \infty)$, there exists
$\bigl(r \mapsto \sum_{n = 0}^\infty a_n^2 r^{2n} \bigr) \in L^{p/2}(\mu)
\setminus \bigcup_{q > p} L^{q/2}(\mu)$.
\endprocl

\rproof
For each $M > 0$, let $q := p + 1/M$ and let $s = s(M) < r_\mu$ be close enough to $r_\mu$
that $\mu(s, r_\mu) \le M^{-pq/(q - p)}$. We can find $N = N(M)$ large enough so that
\[
\int_{s}^{r_\mu} r^{Np} \,d\mu(r) \ge \frac{1}{2} \int_0^{r_\mu} r^{Np} \,d\mu(r)\,.
\]
We also have by the power-mean inequality that
\eqalign{
\Bigl( \int_{s}^{r_\mu} r^{Nq} \,d\mu(r) \Bigr)^{\frac2q}
&\ge
\mu(s, r_\mu)^{\frac2q - \frac2p} \Bigl( \int_{s}^{r_\mu} r^{Np} \,d\mu(r) \Bigr)^{\frac2p}
\\ &\ge
M^2 \Bigl( \int_{s}^{r_\mu} r^{Np} \,d\mu(r) \Bigr)^{\frac2p}.
}

Now, let $n_k := N(2^k)$, and choose $b_{k} > 0$ so that $\left( \int_0^{r_\mu} b_{k}^p r^{n_kp} \,d\mu(r) \right)^{2/p} = 1/2^k$. Then $\sum_{k = 1}^\infty b_{k}^2 r^{2n_k}$ has the desired property: Write
\[
\norm{f}_p := \bigl(\int_0^{r_\mu} |f(r)|^p \,d\mu(r)\bigr)^{1/p}
\,.
\]
If $p \ge 2$, then
\[
\Bigl\|\sum_{k = 1}^\infty b_{k}^2 r^{2n_k}\Bigr\|_{p/2}
\leq
\sum_{k=1}^{\infty}
\norm{b_{k}^2 r^{2n_k}}_{p/2}
=
1\,,
\]
while if $p < 2$, then
\[
\Bigl\|\sum_{k = 1}^\infty b_{n}^2 r^{2n_k}\Bigr\|_{p/2}^{p/2}
\leq
\sum_{k=1}^{\infty}
\norm{b_{k}^2 r^{2n_k}}_{p/2}^{p/2}
<
\infty\,.
\]
At the same time, for each $q > p$, consider any $j$ with $q > p + 1/2^j$. Then
\begin{align*}
\Bigl\|\sum_{k = 1}^\infty b_{k}^2 r^{2n_k}\Bigr\|_{q/2}
&\geq
\bigl\|b_{j}^2 r^{2n_j}\bigr\|_{q/2}
\geq
b_{j}^2 2^{2j} \Bigl( \int_{s(2^j)}^{r_\mu} r^{n_j p} \,d\mu(r) \Bigr)^{\frac2p}
\\ &\geq
b_{j}^2 2^{2j} \Bigl( \frac12 \int_{0}^{r_\mu} r^{n_j p} \,d\mu(r) \Bigr)^{\frac2p}
=
2^{j-2/p}.
\end{align*}
Since this holds for all such $j$, it follows that
$\norm{\sum_{k = 1}^\infty b_{k}^2 r^{2n_k}}_{q/2} = \infty$.
\Qed

\bsection{Proofs}{s.proofs}

In this section, we prove \rref t.quant/ and then indicate the
additional steps needed for the more general \rref t.quant2/.
At the end, we prove \rref t.quant3/.

\begin{proofof} {\rref t.quant/}
Note that the density of $|\zeta_0|$ with respect to
Lebesgue measure on $\R^+$ is $r \mapsto 2 r e^{-r^2}$. It suffices to prove the theorem for $s < r_\mu$, since the case $s = r_\mu$ follows by taking limits.


Suppose that $0 \notin Z(f) \supseteq Z(\sF)$.
Note that $0 \notin Z(\sF)$ a.s. Thus,
for $0 < s < r_\mu$, we have a.s.\
by
the arithmetic mean-geometric mean inequality and Jensen's formula
that
\eqaln{
\norm{f}_{A^p(\mu, s)}^p
&=
\int_0^{s} \int_0^1 |f(r e^{2\pi i \theta})|^p \,d\theta \,d\mu(r)
\\ &\ge
\int_0^{s} \exp \int_0^1 \log |f(r e^{2\pi i \theta})|^p \,d\theta
\,d\mu(r)
\\ &=
\int_0^{s} |f(0)|^p \prod_{z \in Z_s(f)} \max\bigl\{\textfrac{r^p}{|z|^p},
1\bigr\} \,d\mu(r)
\\ &\ge
\int_0^{s} |f(0)|^p \prod_{z \in Z_s(\sF)} \max\bigl\{\textfrac{r^p}{|z|^p},
1\bigr\} \,d\mu(r)
\\ &=
|f(0)/\sF(0)|^p
\int_0^{s} \exp \int_0^1 \log |\sF(r e^{2\pi i \theta})|^p \,d\theta
\,d\mu(r)
\,.
}
Therefore,
\begin{equation}
\frac{|f(0)|}{\norm{f}_{A^p(\mu, s)}}
\le
a_0 |\zeta_0|
\Bigl(\int_0^{s} \exp \int_0^1 \log |\sF(r e^{2\pi i \theta})|^p
\,d\theta \,d\mu(r)\Bigr)^{-1/p}
\,. \label{e.bound-by-F}
\end{equation}

Recall that for each $r$ and $\theta$,
$\> \sF(r e^{2\pi i \theta})$ is a Gaussian random variable with the same
distribution as $\|a^{(r)}\|_2 \zeta_0$, where $a^{(r)}_n := a_n r^n$.
Write
$$
\sG_r(\theta) := \sF(r e^{2\pi i \theta})/\|a^{(r)}\|_2
\,,
$$
so that $\sG_r(\theta)$ is a standard complex Gaussian random variable for each $r$ and $\theta$.
H\"older's inequality and the arithmetic
mean-geometric mean inequality yield
\begin{multline}
\label {e.before}
\Bigl(\int_0^s \exp \int_0^1 \log |\sF(r e^{2\pi i \theta})|^p
\,d\theta \,d\mu(r)\Bigr)^{-1/p} \\
\begin{aligned}
&=
\Bigl(\int_0^s \A r2p \exp \int_0^1 \log |\sG_r(\theta)|^p
\,d\theta \,d\mu(r)\Bigr)^{-1/p}
\\ &\le
\frac{\int_0^s \A r2p \exp \int_0^1 \log |\sG_r(\theta)|^{-1}
\,d\theta \,d\mu(r)
}{
\Bigl(\int_0^s \A r2p \,d\mu(r)\Bigr)^{1+1/p}}
\\ &\le
\frac{\int_0^s \A r2p \int_0^1 |\sG_r(\theta)|^{-1}
\,d\theta \,d\mu(r)
}{
\Bigl(\int_0^s \A r2p \,d\mu(r)\Bigr)^{1+1/p}}\,.
\end{aligned}
\end{multline}
Multiplying both sides by $a_0 |\zeta_0|$ and using \rref e.bound-by-F/, we have
\begin{equation}
\frac{|f(0)|}{\norm{f}_{A^p(\mu, s)}} \le \frac{a_0 |\zeta_0| \cdot \int_0^s \A r2p \int_0^1 |\sG_r(\theta)|^{-1}
\,d\theta \,d\mu(r)
}{
\Bigl(\int_0^s \A r2p \,d\mu(r)\Bigr)^{1+1/p}}\,.
\label {e.bound-by-G}
\end{equation}
Recall that for each $r$ and $\theta$, $\sG_r(\theta)$ and $\zeta_0$ are both standard complex Gaussians, and $(\sG_r(\theta), \zeta_0)$ is jointly Gaussian. By a version of Slepian's lemma due to \rref b.Kahane:slepian/, we have
\begin{equation*}
\Ebig{|\zeta_0|\cdot |\sG_r(\theta)|^{-1}}
\le
\Ebig{|\zeta_0|} \cdot \Ebig{|\sG_r(\theta)|^{-1}}
=
1 \cdot \sqrt\pi
\,.
\end{equation*}
Taking expectations in \rref e.bound-by-G/ and applying the above inequality finishes the proof, except for showing that the maximum on the left-hand side of \rref e.quant/ is achieved and is measurable.

To show these properties,
note first that the maximum is achieved because of a standard normal-families argument (compare \citet[p.~120]{Duren:book}).
Next, for a finite multiset $W$, let $p_W(z) := \prod_{w \in W} (z - w)$ be the monic polynomial whose zeros are $W$ (with multiplicity). For any analytic function $f$ whose zeros include $W$, the function $f/p_W$ is analytic. Therefore,
\begin{multline*}
\max \Bigl\{ \frac{|f(0)| }{ \norm{f}_{A^p(\mu, s)}} \st 0 \not\equiv
f \in A(\mu),\, Z(f) \supset Z_s(\sF) \Bigr\}
\\ =
\max \Bigl\{ \frac{|f(0)p_{Z_s(\sF)}(0)| }{ \norm{fp_{Z_s(\sF)}}_{A^p(\mu, s)}} \st 0 \not\equiv
f \in A(\mu)\Bigr\}\,.
\end{multline*}
Restricting to polynomials $f$ with rational coefficients, we see that this maximum is measurable provided $p_{Z_s(\sF)}$ is measurable.
Now there is a measurable set (of probability 0) where $\limsup |\zeta_n|^{1/n} > 1$; off of this set, $Z_s(\sF)$ is finite and can be determined by
looking at the values of $\sF$ on a fixed, countable, dense set of points, thereby proving the desired measurability.
\end{proofof}

\procl r.noslepian
In fact, \rref t.main/ may be deduced directly from \rref e.before/ without using Slepian's lemma: Simply take expectations of both sides and use the facts that $\Ebig{|\sG_r(\theta)|^{-1}} = \sqrt{\pi}$ and $|\zeta_0| < \infty$ to obtain
\[
\EBig{\Bigl(\int_0^s \exp \int_0^1 \log |\sF(r e^{2\pi i \theta})|^p
\,d\theta \,d\mu(r)\Bigr)^{-1/p}}
\le
\frac{\sqrt\pi
}{
\bigl(\int_0^s \A r2p \,d\mu(r)\bigr)^{1/p}}\,.
\]
As $s \uparrow {r_\mu}$, the right-hand side tends to 0, which already
gives \rref t.main/ via \rref e.bound-by-F/.
\endprocl


\begin{proofof} {\rref t.quant2/}
We may assume that $a_0 \ne 0$.
Suppose that $0 \notin Z(f) \supseteq Z(\sF)$.
As before, we have for $0 < s < r_\mu$
\begin{multline*}
\Bigl(\frac{|f(0)|}{\norm{f}_{A^{p/N}(\mu, s)}}\Bigr)^{1/N}
\le |a_0^N \zeta_0^N + b(0)|^{1/N} \cdot \\
\Bigl(\int_0^{s} \exp \int_0^1 \log |\sF(r e^{2\pi i \theta})^N + b(r
e^{2\pi i \theta})|^{p/N}
\,d\theta \,d\mu(r)\Bigr)^{-1/p}
\,.
\end{multline*}
Write
$$
\sH_r(\theta) := |\sF(r e^{2\pi i \theta})^N + b(r e^{2\pi i \theta})|/\A r2N
\,.
$$
In the same way as before, we obtain
\begin{multline}
\label {e.before2}
\Bigl(\int_0^s \exp \int_0^1 \log |\sF(r e^{2\pi i \theta})^N +
b(r e^{2\pi i \theta})|^{p/N}
\,d\theta \,d\mu(r)\Bigr)^{-1/p} \\
\begin{aligned}
&=
\Bigl(\int_0^s \A r2p \exp \int_0^1 \log \sH_r(\theta)^{p/N}
\,d\theta \,d\mu(r)\Bigr)^{-1/p}
\\ &\le
\frac{\int_0^s \A r2p \exp \int_0^1 \log \sH_r(\theta)^{-1/N}
\,d\theta \,d\mu(r)
}{
\Bigl(\int_0^s \A r2p \,d\mu(r)\Bigr)^{1+1/p}}
\\ &\le
\frac{\int_0^s \A r2p \int_0^1 \sH_r(\theta)^{-1/N}
\,d\theta \,d\mu(r)
}{
\Bigl(\int_0^s \A r2p \,d\mu(r)\Bigr)^{1+1/p}}\,.
\end{aligned}
\end{multline}
We have for any $\beta$ that
$$
\Ebig{\sH_r(\theta)^{-\beta}}
=
\Ebig{|\zeta_0^N + b(r e^{2\pi i \theta})/\A r2N|^{-\beta}}
\,.
$$
Now $\zeta_0^N$ has
density $\rho \colon z \mapsto c |z|^{-2(N-1)/N} e^{-|z|^2}$ (with respect to area measure $\lambda_2$, for some constant $c$) that is decreasing in $|z|$. Therefore, given any $\alpha \in \C$, the rearrangement inequality of Hardy and Littlewood yields
\[
\Pbig{|\zeta_0^N| < r} = \int_{|z| < r} \rho(z) \,d\lambda_2(z) \ge \int_{|z - \alpha| < r} \rho(z) \,d\lambda_2(z) = \Pbig{|\zeta_0^N - \alpha| < r}\,,
\]
which is to say that $|\zeta_0^N|$ is stochastically dominated by $|\zeta_0^N - \alpha|$.
Thus, for all $0 < \beta < 2/N$,
\[
\Ebig{|\zeta_0^N - \alpha|^{-\beta}}
\le
\Ebig{|\zeta_0^N|^{-\beta}}
=
\int_0^\infty \frac{e^{-t} }{ t^{\beta N/2}} dt
=
\Gamma(1 - \beta N/2)
\,.
\]
Therefore,
\rlabel e.PMineq
{
\Ebig{\sH_r(\theta)^{-\beta}}
\le
\Gamma(1 - \beta N/2)
\,.
}
Multiply both sides of the
inequality \rref e.before2/ by $|a_0^N \zeta_0^N + b(0)|^{1/N}$
and use H\"older's inequality
to bound the resulting expectation:
\begin{multline*}
\Ebig{|a_0^N \zeta_0^N + b(0)|^{1/N}\cdot
\sH_r(\theta)^{-1/N}} \\
\begin{aligned}
&\le
\Ebig{|a_0^N \zeta_0^N + b(0)|^{4}}^{\frac1{4N}}
\Ebig{\sH_r(\theta)^{-4/(4N-1)}}^{\frac{4N - 1}{4N}}
\\ &\le
\Bigl(a_0^{4N} (2N)! + 4|b(0)|^2 a_0^{2N} N! + |b(0)|^4\Bigr)^{1/4N}
\Gamma\Bigl(\frac{2N-1}{4N-1}\Bigr)^{\frac{4N - 1}{4N}}
\,,
\end{aligned}
\end{multline*}
where in the last inequality, we used \rref e.PMineq/ with $\beta :=
4/(4N-1)$.
\end{proofof}

\begin{proofof}{\rref t.quant3/}
We established \rref e.quant3/ during the proof of \rref t.quant/, where we rely on
\rref e.tonelli/ and the fact that $\Ebig{|\zeta_0|^p} = \Gamma(1+p/2)$ for an equivalent expression on the right-hand side.
That (ii) implies (iii) follows from the arithmetic mean-geometric mean inequality.
That (iii) implies (i) and (i) implies (iv) are obvious. That (iv) implies (ii) follows from \rref e.quant3/ with $s = R$.
\end{proofof}

\appendix
\bsection{Fernique's Theorem}{s.fernique}

We present here a general version of Fernique's theorem, not only for use in deriving the background in \autoref{sec:shap-app}, but also for comparison with our \rref t.quant3/.

\procl t.fernique
Let $V$ be a separable topological vector space.
Let $\phi \colon V \to [0, \infty]$ be Borel measurable,
$\>c \in [1, \infty)$, and $c_1, c_2 \in (1, \infty)$ satisfy
for all $x, y \in V$ that
$\phi(-x) = \phi(x)$,
$\>c_2 \phi(x) \le \phi(\sqrt2\,x) \le c_1 \phi(x)$,
and $\phi(x + y) \le c\bigl(\phi(x) + \phi(y)\bigr)$.
Let $X$ be a random variable with values in $V$ such
that if $Y$ has the same distribution as $X$ and is independent of $X$,
then $\bigl(\phi(X), \phi(Y)\bigr)$ has the same joint distribution as
$\bigl(\phi\bigl((X-Y)/\sqrt2\bigr),
\phi\bigl((X+Y)/\sqrt2\bigr)\bigr)$.
If\/ $\Pbig{\phi(X) < \infty} = 1$,
then there are some $\alpha, \beta > 0$ so that
$\Ebig{e^{\alpha \phi(X)^\beta}} < \infty$.
\endprocl

\rproof
Suppose that $\phi\bigl((x-y)/\sqrt2\bigr) \le \tau$ and
$\phi\bigl((x+y)/\sqrt2\bigr) > t$.
Then $\phi(x-y) \le c_1 \tau$ and $\phi(x+y) > c_2 t$.
Also, $\phi(2y) = \phi(\sqrt2^2y) \le c_1^2 \phi(y)$, whence
$\phi(x + y) \le c \phi(x - y) + c c_1^2 \phi(y)$.
Therefore $\phi(y) > (c_2 t - c c_1 \tau)/(cc_1^2)$.
Symmetry gives the same lower bound on $\phi(x)$.
It follows that
\begin{align}
\label {e.rotate}
\begin{split}
\Pbig{\phi(X) \le \tau} &\Pbig{\phi(Y) > t}
\\ &=
\Pbig{\phi\bigl((X-Y)/\sqrt2\bigr) \le \tau,\,
\phi\bigl((X+Y)/\sqrt2\bigr) > t}
\\ &\le
\Pbig{\phi(X) > (c_2 t - c c_1 \tau)/(cc_1^2)}^2
\,.
\end{split}
\end{align}
Choose $\tau < \infty$ so that $\Pbig{\phi(X) \le \tau} \ge e/(1+e)$.
Define recursively $t_0 := (c c_1/c_2)\tau \ge \tau$ and $t_{n+1} :=
(cc_1^2/c_2) t_n + t_0$;
thus,
\[
t_n = \frac{(cc_1^2/c_2)^{n+1} - 1 }{ cc_1^2/c_2 - 1} t_0
<
c_3 (cc_1^2/c_2)^n t_0
\]
for some constant $c_3 < \infty$.
The display \rref e.rotate/ yields
\[
\Pbig{\phi(X) > t_{n+1}}
=
\Pbig{\phi(Y) > t_{n+1}}
\le
\frac{1+e }{ e} \Pbig{\phi(X) > t_n}^2
\,,
\]
whence if we write $y_n := \frac{1+e }{ e} \Pbig{\phi(X) > t_n}$, then
$y_{n+1} \le y_n^2$, and so $y_n \le y_0^{2^n} \le e^{-2^n}$.
Therefore,
\[
\Pbig{\phi(X) > c_3 (cc_1^2/c_2)^n t_0}
\le
e^{-2^n}
=
e^{-(cc_1^2/c_2)^{\beta n}}
\,,
\]
where $\beta := \log 2/\log(cc_1^2/c_2) > 0$.
This means that
\[
\Pbig{\phi(X) > t}
\le
e^{-c_4 t^\beta}
\]
for some $c_4 > 0$ and all $t \ge t_0$.
With $\alpha := c_4/2$, the conclusion may be verified via integration
by parts.
\Qed

\nocite{*}
\bibliographystyle{cdraifplain}

\pdfbookmark{References}{refs}

\bibliography{zeros}

\end{document}